\documentclass[12pt]{amsart}
\usepackage{amssymb}
\usepackage[dvips]{graphics}
\ifx\mathrm\undefined\let\mathrm\rm\fi
\ifx\mathbf\undefined\let\mathbf\bf\fi
\ifx\mathfrak\undefined\let\mathfrak\frak\fi
\ifx\mathcal\undefined\let\mathcal\cal\fi
\ifx\mathbb\undefined\let\mathbb\Bbb\fi
\ifx\emph\undefined\let\emph\it\fi
 at9.98pt

\newcommand{\g}{{{\mathfrak g}\,}}

\newcommand{\n}{{{\mathfrak n}}}

\newcommand{\h}{{{\mathfrak h\,}}}

\newcommand{\Z}{{\mathbb Z}}

\newcommand{\C}{{\mathbb C}}

\newcommand{\Ref}[1]{{(\ref{#1})}}

\newcommand{\la}{\lambda}

\newcommand{\dontprint}[1]
{\relax}

\newtheorem%
{thm}{Theorem}[section]
\newtheorem%
{proposition}[thm]{Proposition}
\newtheorem%
{lemma}[thm]{Lemma}
\newtheorem%
{lemmadef}[thm]{Lemma-Definition}
\newtheorem%
{corollary}[thm]{Corollary}
\newtheorem%
{conjecture}[thm]{Conjecture}

\newcommand{\bea}{\begin{eqnarray*}}
\newcommand{\eea}{\end{eqnarray*}}
\newcommand{\bean}{\begin{eqnarray}}
\newcommand{\eean}{\end{eqnarray}}

\newcommand{\nc}{\newcommand}
\nc{\on}{\operatorname}
\nc{\al}{\alpha}
\nc{\ri}{\rangle}
\nc{\lef}{\langle}
\nc{\W}{{\mathcal W}}
\nc{\La}{\Lambda}
\nc{\ep}{\epsilon}
\nc{\Om}{\Omega}
\newcommand{\be}{\begin{displaymath}}
\newcommand{\ee}{\end{displaymath}}
\newcommand{\bs}{\boldsymbol}
\nc{\PCr}{{ \Bbb P  (\C[x])^r   }}
\newtheorem{theorem}{Theorem}[section]

\newtheorem{cor}[theorem]{Corollary}

\newcommand{\id}{{{\rm id}}}
\newcommand{\om}{{\rm{Om}_{D_{\bs y}}}}

\newcommand{\mc}{\mathcal}

\newcommand{\s}{{\rm Sing\,}}

\newcommand{\End}{{\rm End\,}}

\newcommand{\V}{ V_{\bs \La}}

\newcommand{\Ll}{\La - \al(\bs l)}

\newcommand{\mult}{{\rm mult}\,}

\begin{document}

\title[The Norm of a Bethe Vector and the Hessian  of the Master Function]
{Norm of a Bethe Vector and the Hessian \\
of the Master Function}

\author[{}]{ Evgeny Mukhin  
${}^{*,1}$ 
\and Alexander Varchenko${}^{**,2}$}
\thanks{${}^1$ Supported in part by NSF grant DMS-0140460}

\thanks{${}^2$ Supported in part by NSF grant DMS-0244579}

\begin{abstract}
We show that the Bethe vectors are non-zero vectors in the $sl_{r+1}$ Gaudin model.
Namely, we show that the norm of a Bethe vector is equal to the Hessian
of the corresponding master function at the corresponding non-degenerate
critical point. 
This result is a byproduct of functorial properties  of Bethe vectors 
studied in this paper. 

As other byproducts of functoriality we
show that the Bethe vectors form a basis in the tensor product of
several copies of first and last fundamental $sl_{r+1}$ modules and we
show transversality of some  Schubert cycles in the Grassmannian
of $r+1$-dimensional planes in the space $\C_d[x]$ of polynomials of one variable of
degree not greater than $d$.

\end{abstract}

\maketitle 
\medskip \centerline{\it ${}^*$
Department of Mathematical Sciences,}
\centerline{\it Indiana University
Purdue University Indianapolis,}
\centerline{\it 402 North Blackford St., Indianapolis,
IN 46202-3216, USA}
 \medskip
\centerline{\it ${}^{**}$Department of Mathematics, University of
  North Carolina at Chapel Hill,} \centerline{\it Chapel Hill, NC
  27599-3250, USA} \medskip

\centerline{February, 2004}

\section{Introduction}

The Bethe ansatz is a large collection of methods in the theory of
quantum integrable models to calculate the spectrum and eigenvectors
for a certain commutative sub-algebra of observables for an integrable
model. Elements of the sub-algebra are called hamiltonians, or integrals of motion,
or conservation laws of the model.  The bibliography on
the Bethe ansatz method is enormous, see for example
\cite{BIK, Fa, FT}. 

In the theory of the Bethe ansatz one assigns the Bethe ansatz equations to
an integrable model. Then a solution of the Bethe ansatz equations gives an
eigenvector of commuting hamiltonians of the model. The general
conjecture is that the constructed vectors form a basis in the space of 
states of the model.

The simplest and interesting example is the Gaudin model associated
with a complex simple Lie algebra $\g$, see \cite{B, BF, F1, FFR, G, MV1, RV, ScV, V2}.
One considers highest weight $\g$-modules
$V_{\La_1}, \dots , V_{\La_n}$ and their tensor product $\V$. One fixes
a point $z = (z_1, \dots , z_n) \in \C^n$ with distinct coordinates
and defines linear operators
$K_1(z), \dots , K_n(z)$ on $V_{\bs \La}$ by the formula
\bea
K_i(z)\ = \ \sum_{j \neq i}\ \frac{\Omega^{(i,j)}}{z_i - z_j} , 
\qquad i = 1, \dots , n .
\eea
Here $\Omega^{(i,j)}$ is the Casimir operator acting in the $i$-th and $j$-th factors 
of the tensor product. The operators are called  the Gaudin hamiltonians
of the Gaudin model associated with $\V$. The hamiltonians commute.

The common eigenvectors of the Gaudin hamiltonians 
are constructed by the Bethe ansatz method. Namely, one assigns to the model
a scalar function $\Phi (t,z)$ of new auxiliary variables $t$ and a
$\V$-valued function $\omega(t,z)$ such that $\omega(t^0,z)$ is  an eigenvector
of the hamiltonians if $t^0$ is a critical point of  $\Phi$.
The functions $\Phi$ and $\omega$ were introduced in \cite{SV} to construct 
hypergeometric solutions of the KZ equations. 
The function $\Phi$ is called the { master function} and the function
$\omega$ is called { the universal weight function}.

The first question is
if the Bethe eigenvector $\omega(t^0,z)$ is non-zero. In this paper we
show that for the $sl_{r+1}$ Gaudin model the Bethe vector
is non-zero if $t^0$ is a non-degenerate critical point of the master function $\Phi$.
To show that we prove the following identity:
\bean\label{first}
S ( \omega (t^0,z), \omega (t^0,z) )\ =\
\text{Hess}_t\ \text{log}\ \Phi (t^0, z) \ .
\eean
Here $S$ is the tensor Shapovalov form on the tensor product $\V$ and
the right hand side of the formula is the Hessian at $t^0$ of the function
${\rm log}\ \Phi$. This formula for $sl_2$ Gaudin models was proved in 
\cite{V2}, see also \cite{RV, Ko, R, TV, MV1}.

\bigskip

In this paper we prove 
the Bethe ansatz conjecture for tensor products of several copies 
of first and last fundamental $sl_{r+1}$-modules. Namely,  if
$V_{\La_1}, \dots , V_{\La_n}$ are  $sl_{r+1}$-modules, 
each of which is either the first or last fundamental $sl_{r+1}$-module,
then we show that for generic $z$ the Bethe vectors form an eigenbasis 
of the Gaudin hamiltonians in the tensor product $\V$. 
Note that $sl_3$ has only two fundamental modules: the first and last.

We also prove the Bethe ansatz conjecture for tensor products of
several copies of arbitrary fundamental representations of $sl_4$.

\bigskip

The formulated results are based on functorial properties of the
master function and the universal weight function studied in this paper.
Namely we study the behavior of $\Phi$ and $\omega$ when some of coordinates
of $z$ tend to the same limit. That corresponds to the situation in which
the number of factors in the tensor product $\V$ 
becomes smaller while the factors become bigger.
It turns out that under this limit the Bethe vectors behave in a reasonable way.
That reasonable behavior  allows us to establish some general properties of Bethe vectors
under the condition that   those properties hold for some model examples.
The properties for the model examples can be checked by direct calculations. 

The ideas of that type were exploited earlier in \cite{RV}.

\bigskip

The paper is organized as follows. 
Section 2 contains the definition of the master  and universal weight functions.
We prove there that the universal function is well defined on critical 
points of the master function. In Section 3 we collect information on iterated
singular vectors in tensor products of representations.
The functorial properties of the master and universal weight functions
are studied in Section 4. Preliminary information on Bethe vectors and their Shapovalov 
norms is collected in Section 5. In Section 6 we prove Theorem 
\ref{sl3 fund weights} that the Bethe vectors form a basis
in the tensor product of several copies of
first and last fundamental  $sl_{r+1}$-modules for generic $z$.
In Section 7 we prove formula \Ref{first} using Theorem
\ref{sl3 fund weights}.
In Section 8 as a corollary of Theorem \ref{sl3 fund weights}
we show transversality of some Schubert cycles in the Grassmannian
of $r+1$-dimensional planes in the space $\C_d[x]$ of polynomials of one variable of
degree not greater than $d$.

\section{Bethe Vectors}\label{vectors}

\subsection{The Gaudin model}\label{model}
Let $\g$ be a simple Lie algebra over $\C$ with Cartan matrix
$A=(a_{i,j})_{i,j=1}^r$. Let $D=\on{diag}\{d_1,\dots,d_r\}$ be the diagonal
matrix with positive relatively prime
integers $d_i$ such that $B=DA$ is symmetric.

Let $\h \subset \g$  be the Cartan sub-algebra.
Fix simple roots $\al_1, \dots , \al_r$ in $\h^*$
and an invariant bilinear form $( , )$ on $\g$ such that
$(\al_i , \al_j) = d_i a_{i,j}$. Let $H_1, \dots , H_r \, \in \h$ be the corresponding
coroots, $\langle\la , H_i\rangle = 2 (\la,\al_i)/ (\al_i,\al_i)$
 for $\la\in\h^*$.
In particular, $\langle \al_j , H_i \rangle = a_{i,j}$.
Let $w_1, \dots , w_r \in \h^*$
be the fundamental weights, $ \langle w_i , H_j \rangle = \delta_{i,j}$.

Let $E_1, \dots , E_r\, \in \n_+,\ H_1, \dots , H_r\, \in \h,\
F_1, \dots , F_r\, \in \n_- $ be the Chevalley generators of $\g$,
\bea
[E_i, F_j] & = & \delta_{i,j} \,H_i ,
\qquad i, j = 1, \dots r ,
\\
{}[ h , h'] & = & 0 ,
\qquad 
\phantom{aaaaa}
h, h' \in \h ,
\\
{}[ h, E_i] &=& \langle \alpha_i, h \rangle\, e_i ,
\qquad  
h \in \h, \ i = 1, \dots r ,
\\
{}[ h, F_i] &=& - \langle \alpha_i, h \rangle\, F_i ,
\qquad  h\in \h, \ i = 1, \dots r ,
\eea
and
\bea
(\mathrm{ad}\,{} E_i)^{1-a_{i,j}}\,E_j = 0 ,
\qquad
\phantom{aaaaaaa}
(\mathrm{ad}\, {} F_i)^{1-a_{i,j}}\,F_j = 0 ,
\eea
for all $i\neq j$.

Let $(x_i)_{i\in I}$ be an orthonormal basis in $\g$, \
$\Omega =  \sum_{i\in I} x_i\otimes x_i\ \in \g \otimes \g$
the Casimir element. We have
\bean\label{casinir property}
[x\otimes 1 + 1 \otimes x, \ \Omega]\ =\ 0
\eean
in $U(\g) \otimes U(\g)$ for any $x \in \g$. Here $U(\g)$ is the universal 
enveloping algebra of $\g$. 

For a $\g$-module $V$ and $\mu \in \h^*$
denote by $V[\mu]$ the weight subspace of $V$ of weight $\mu$ and by
$\s V[\mu]$ the subspace of singular vectors of weight $\mu$,
\bea
\s V[\mu]\ =\ \{ \ v \in V\ |\ \n_+v = 0, \ hv = \langle \mu, h \rangle v \ \} \ .
\eea

Let $n$ be a positive integer and  $\bs \La = (\La_1, \dots , \La_n)$,
$\La_i \in \h^*$, a set of weights.
For $\mu \in \h^*$ let $V_{\mu}$ 
be the irreducible $\g$-module with highest weight $\mu$.
Denote by $V_{\bs \La}$ the tensor product 
$V_{\La_1} \otimes \dots \otimes V_{\La_n}$.

If $X \in \End\,(V_{\La_i})$, then we denote by $X^{(i)} \in \End (V_{\bs \La})$ the operator
$ \cdots \otimes \id \otimes X \otimes \id \otimes \cdots$ acting non-trivially on the
$i$-th factor of the tensor product. If $X = \sum_k X_k \otimes Y_k \in 
\End (V_{\La_i} \otimes V_{\La_j})$, then we set 
$X^{(i,j)} = \sum_k X^{(i)}_k \otimes Y^{(j)}_k\ \in \End (\V)$.


Let $z = (z_1, \dots , z_n)$ be a point in $\C^n$ with distinct coordinates.
Introduce linear operators $K_1(z), \dots , K_n(z)$ on $V_{\bs \La}$ by the formula
\bea
K_i(z)\ = \ \sum_{j \neq i}\ \frac{\Omega^{(i,j)}}{z_i - z_j}\ , 
\qquad i = 1, \dots , n .
\eea
The operators 
are called {\it the Gaudin hamiltonians} of the Gaudin model associated with
$\V$. One can check directly
that the  hamiltonians commute, $ [ H_i(z), H_j(z) ] = 0$ for all $i, j$.

The main problem for the Gaudin model is to diagonalize simultaneously
the  hamiltonians, see \cite{B, BF, F1, FFR, G, MV1, RV, ScV, V2}.

One can check that the  hamiltonians commute with the action of $\g$ on
$V_{\bs \La}$, $[ H_i(z), x ] = 0$ for all $i$ and $x \in \g$. Therefore it is enough to
diagonalize the  hamiltonians on the subspaces of singular vectors
$\s V_{\bs \La}[\mu] \subset \V$. 

The eigenvectors of the Gaudin hamiltonians are constructed by the Bethe ansatz method.
We  remind the construction in the next section.

\subsection{Master functions, critical points, and the universal weight function}
\label{master sec} 
Fix  a collection of weights $\bs\La = (\La_1, \dots , \La_n)$, $\La_i\in  \h^*$,
and a collection of non-negative integers  $\bs l = (l_1, \dots , l_r)$. 
Denote  $l = l_1 + \dots + l_r$, $\La = \La_1 + \dots + \La_n$, and 
$\al(\bs l) = l_1\al_1 + \dots + l_r\al_r$.

Let $c$ be the unique non-decreasing function from 
$\{1, \ldots , l\}$ to $\{1, \ldots , r\}$, such that 
$\# c^{-1}( i) = l_i$ for $i = 1, \dots , r$. 
The {\it master function} $\Phi(t, z, \bs \La, \bs l)$
is defined by the formula
\bea\label{master}
\Phi(t,  z, \bs \La, \bs l) = \prod_{1\leq i < j\leq n}
  (z_i - z_j)^{(\La_i, \La_j)}
 \prod_{i=1}^l \prod_{s=1}^n
  (t_i - z_s)^{-(\al_{c(i)}, \La_s)}
  \prod_{1 \leq i < j \leq l}  (t_i - t_j)^{(\al_{c(i)},\al_{c(j)})} ,
\eea
see  \cite{SV}.
The function $\Phi$ is a function of complex variables 
$t = ( t_1, \dots , t_l)$, $z = (z_1, \dots , z_n)$, weights $\bs \La$,
and discrete parameters $\bs l$.
The main variables are $t$, the other variables will be
considered as parameters.

For given $z,  \bs \La, \bs l$, a point $t$ with complex coordinates is called {\it a critical
point} of the master function
if the system of algebraic equations is satisfied,
\bean\label{Bethe eqn}
- \sum_{s=1}^n \frac{(\alpha_{c(i)}, \La_s)}{t_i - z_s}\ +\
\sum_{j,\ j\neq i} \frac{(\alpha_{c(i)}, \alpha_{c(j)})}{ t_i - t_j}
= 0, 
\qquad
i = 1, \dots , l .
\eean
In other words,  $t$ is a critical point if
\be
\left(\Phi^{-1}\frac{\partial \Phi }{\partial t_i}\right) ( t )\ =\ 0, \qquad
{\rm for} \ {} \ i = 1, \dots , l\ .
\ee
By definition, if $t=(t_1, \dots , t_l)$ is a critical point and 
$(\al_{c(i)}, \al_{c(j)})\neq 0$ for some $i, j$, then $t_i \neq t_j$. Also if 
$(\al_{c(i)}, \La_s)\neq 0$ for some $i, s$, then $t_i \neq z_s$.

Let $\Sigma_l$ be the permutation group of the set $\{1, \dots , l\}$.
Denote by $\bs \Sigma_{\bs l} \subset \Sigma_l$ the subgroup of all permutations preserving
the level sets of the function $c$. The subgroup $\bs\Sigma_{\bs l}$ is isomorphic to 
$\Sigma_{l_1}\times \dots \times \Sigma_{l_r}$ and 
acts on $\C^l$ permuting coordinates of $t$. The
action of the subgroup $\bs\Sigma_{\bs l}$ preserves the  critical set
of the master function. All orbits of $\bs\Sigma_{\bs l}$ on the critical set
have the same cardinality $l_1! \cdots l_r!$\ .

\bigskip

Consider highest weight irreducible $\g$-modules
$V_{\La_1}, \dots , V_{\La_n}$, the tensor product $\V$,
and its weight subspace $\V[\La - \al(\bs l)]$. Fix a highest weight vector
$v_{\La_i}$ in $V_{\La_i}$ for for all $i$.

We construct a rational map
\bea
\omega \ :\ \C^l \times \C^n\ \to \V[\Ll]
\eea
called {\it the universal weight function}.

Let $P(\bs l,n)$ be the set of sequences 
$I\ = \ (i_1^1, \dots , i^1_{j_1};\ \dots ;\ i^n_1, \dots , i^n_{j_n})$ of integers in 
$\{1, \dots , r\}$ 
such that for all
$i = 1, \dots ,  r$, the integer $i$ appears in $I$ precisely $l_i$ times. 
For $I \in P(\bs l, n)$, and a permutation $\sigma \in \Sigma_l$, 
set $\sigma_1(i) = \sigma(i)$ for $i = 1, \dots , j_1$,
and $\sigma_s(i) = \sigma(j_1+\cdots +j_{s-1}+i)$ for $s = 2, \dots , n$ and 
$ i = 1, \dots , j_s$. Define 
\bea
\Sigma(I)\ {} = \ {}
\{\ \sigma \in \Sigma_l\ {} |\ {} c(\sigma_s(j)) = i_s^j 
\ {} \text{for} \ {} s = 1, \dots , n \ {} \text{and} \ {} j = 1, \dots j_s\ \}\ .
\eea 

To every $I \in P(\bs l, n)$ we associate a vector 
\bea
F_Iv\ =\ F_{i_1^1} \dots F_{i_{j_1}^1}v_{\La_1} \otimes \cdots
\otimes F_{i_1^n} \dots F_{i_{j_n}^n}v_{\La_n}
\eea
in $\V[\Ll]$, and  rational functions
\bea
\omega_{I,\sigma} \ =\ \omega_{\sigma_1(1),\ldots,\sigma_1(j_1)}(z_1)\
\cdots\
\omega_{\sigma_n(1),\ldots,\sigma_n(j_n)}(z_n) ,
\eea
labeled by $\sigma\in \Sigma(I)$,
where 
\bea
\omega_{i_1,\ldots, i_j}(z_s) \ =\ \frac 1 {(t_{i_1}-t_{i_2}) \cdots
(t_{i_{j-1}}-t_{i_j}) (t_{i_j}-z_s)} .
\eea
We set
\bean\label{bethe vector}
\omega(z,t) \ =\  \sum_{I\in P(\bs l,n)}\ \sum_{\sigma\in \Sigma (I)}\ 
\omega_{I,\sigma}\ F_I v\ .
\eean

\noindent
{\bf Examples.} If $\bs l =  (1, 1, 0, \dots , 0)$, then
\bea
\omega(t,z) = 
\frac 1{(t_1-t_2)(t_2-z_1)}F_1F_2v_{\La_1}\otimes v_{\La_2} +
\frac 1{(t_2-t_1)(t_1-z_1)}F_2F_1v_{\La_1}\otimes v_{\La_2} 
\\
+
\frac 1{(t_1-z_1)(t_2-z_2)}F_1v_{\La_1}\otimes F_2v_{\La_2} +
\frac 1{(t_2-z_1)(t_1-z_2)}F_2v_{\La_1}\otimes F_1v_{\La_2} 
\\
+
\frac 1{(t_1-t_2)(t_2-z_2)}v_{\La_1}\otimes F_1F_2v_{\La_2} +
\frac 1{(t_2-t_1)(t_1-z_2)}v_{\La_1}\otimes F_2F_1v_{\La_2} .
\eea
 If $\bs l =  (2, 0, \dots , 0)$, then
\bea
\omega(t,z) 
&=& 
(\frac 1{(t_1-t_2)(t_2-z_1)}+\frac 1{(t_2-t_1)(t_1-z_1)})
\ {}F_1^2v_{\La_1}\otimes v_{\La_2}
\\
& +&
(\frac 1{(t_1-z_1)(t_2-z_2)}+\frac 1{(t_2-z_1)(t_1-z_2)})
\ {} F_1v_{\La_1}\otimes F_1v_{\La_2}
\\
& +&
(\frac 1{(t_1-t_2)(t_2-z_2)}+\frac 1{(t_2-t_1)(t_1-z_2)})
\ {} \ {}\ v_{\La_1}\otimes F_1^2v_{\La_2} .
\eea

The universal weight function  was introduced in \cite{SV}  to solve 
the KZ equations, see \cite{SV, FSV2, FMTV}. The hypergeometric solutions to the KZ equations
with values in $\s  \V[\Ll]$ have the form
\bea
I(z)\ =\ \int_{\gamma(z)} \ \Phi(t,  z, \bs \La, \bs l)^{1/\kappa} \ \omega(t,z)\ dt .
\eea
The values of the universal function are called {\it the Bethe vectors}, 
see \cite{RV, V2, FFR}.

\begin{lemma}\label{well def}
Assume that $z \in \C^n$ has distinct coordinates.
Assume that  $t \in \C^l$ is a critical point of the master function
$\Phi(\, .\, ,  z, \bs \La, \bs l)$. Then the  vector
$\omega(t,z) \in \V[\Ll]$ is well defined. 
\end{lemma}

\begin{proof} The rational function $\omega$ of $t$ and $z$ may have poles 
at hyperplanes given by equations of the form $t_i-t_j=0$ and $t_i-z_s=0$.
All of the poles are of  first order.
We need to prove two facts:
\begin{enumerate}
\item If $(\al_{c(i)}, \al_{c(j)}) = 0$ for some $i$ and $j$, then $w$ does not have
a pole at the hyperplane $t_i-t_j=0$.
\item If $(\al_{c(i)}, \La_s) = 0$ for some $i$ and $s$, 
then $w$ does not have a pole at the hyperplane $t_i-z_s=0$.
\end{enumerate}
Assume that $(\al_{c(i)}, \al_{c(j)}) = 0$ for some $i$ and $j$. 
From formulas for $\omega_{I,\sigma}$ it follows that 
the residue of $\omega$ at  $t_i-t_j=0$ belongs to the span of the vectors
in $\V$ having the form
\bea
 F_{i_1^1} \dots F_{i_{j_1}^1}v_{\La_1} \otimes \cdots \otimes F_{i_{1}^s} \dots 
(F_{c(i)}F_{c(j)}-F_{c(j)}F_{c(i)})\dots 
F_{i_{j_s}^s} v_{\La_s} \otimes  \cdots
\otimes F_{i_1^n} \dots F_{i_{j_n}^n}v_{\La_n} .
\eea
But the element  $F_{c(i)}F_{c(j)}-F_{c(j)}F_{c(i)}$ acts by zero on $\V$. 
Hence $\omega $ is regular at $t_i-t_j=0$.

Assume that $(\al_{c(i)}, \La_s) = 0$ for some $i$ and $s$. 
From formulas for $\omega_{I,\sigma}$ it follows that 
the residue of $\omega$ at  $t_i-z_s=0$ belongs 
to the span of monomials 
\bea
F_Iv \ = \ \cdots \ \otimes\ F_{i_{1}^s} \dots 
F_{i_{j_s}^s} v_{\La_s}\ \otimes \ \cdots
\eea
such that $F_{i_{j_s}^s}=F_{c(i)}$. But $F_{c(i)}v_{\La_s} = 0$
in the irreducible $\g$-module $V_{\La_s}$. 
Hence $\omega $ is regular at $t_i-z_s=0$.
\end{proof}

\begin{theorem}[\cite{RV}]\label{cr bethe}
Assume that $z \in \C^n$ has distinct coordinates.
Assume that  $t \in \C^l$ is a critical point of the master function
$\Phi(\, .\, ,  z, \bs \La, \bs l)$. Then 
the vector $\omega(t,z)$ belongs to $\s\V[\Ll]$ and
 is an eigenvector of the Gaudin hamiltonians $K_1(z), \dots , K_n(z)$.
\end{theorem}

This theorem was proved in \cite{RV} using the  
quasi-classical asymptotics of the hypergeometric solutions
of the KZ equations. The theorem also follows 
directly from Theorem  6.16.2 in \cite{SV},
cf. Theorem 7.2.5 in \cite{SV},\ see also Theorem 4.2.2 in \cite{FSV2}.

\section{The Shapovalov Form and Iterated Singular Vectors}

\subsection{The Shapovalov Form}
Define the anti-involution $\tau : \g \to \g $ sending
$E_1, \dots , E_r, \linebreak H_1, \dots , H_r, \
F_1, \dots , F_r$ to
$F_1, \dots , F_r, \ H_1, \dots , H_r, \
E_1, \dots , E_r$, respectively.

Let $W$ be a highest weight $\g$-module   with highest weight vector $w$.
{\it The Shapovalov form}  on $W$ is the unique
symmetric bilinear form $S$
defined by the conditions:
\bea
S(w, w) = 1 ,
\qquad
S(xu, v) = S(u, \tau(x)v)
\eea
for all $u,v \in W$ and $x \in \g$.

Let $V_{\La_1}, \dots , V_{\La_n}$ be irreducible
highest weight modules and
$\V$ their tensor product.
Let $v_{\La_i} \in V_{\La_i}$ be 
a highest weight vector and 
$S_i$ the corresponding Shapovalov form on $V_{\La_i}$. 
 Define 
a symmetric bilinear form on $\V$  by the formula
\bean\label{shap}
S\ =\ S_1 \otimes \cdots \otimes S_n .
\eean
The form $S$ will be called {\it the tensor Shapovalov form on $\V$}.

\begin{lemma}[\cite{RV}]\label{chap and gaudin}
The Gaudin hamiltonians $K_1(z), \dots , K_n(z)$ are
symmetric with respect to $S$,
$S(K_i(z)u, v) = S(u, K_i(z)v)$ for all $i, z, u, v$.
\end{lemma}

\subsection{Iterated singular vectors}\label{iterated}
Let $n_1, \dots , n_k$ be positive integers.
For $p = 0, 1, \dots , k$ fix  a collection of non-negative integers
$\bs l^p = (l_1^p, \dots , l_r^p)$.
Set $\bs l = \bs l^0 + \bs l^1 + \dots + \bs l^k$, \
  $\al(\bs l^p) = l^p_1\al_1 + \dots + l^p_r\al_r$, \
$n = n_1 + \dots + n_k$, \
$l^p = l^p_1 + \dots + l^p_r$, \
$l = l^0 + l^1 + \dots + l^k$.
For $j = 1, \dots , r$, set
$l_j = l^0_j + l^1_j + \dots + l^k_j$. We have $l = l_1 + \dots + l_r$.

For $p = 1, \dots , k$ fix  a collection of weights 
$\bs \La^p = (\La_1^p, \La_2^p, \dots , \La^p_{n_p}),\ \La_i^p \in \h^*$.
Denote by $\bs \La$ the collection of $n$ weights $\La^p_i$, $p = 1, \dots , k, \
i = 1, \dots , n_p$.
Set $\La^p = \La_1^p + \dots + \La_{n_p}^p$, 
$\La = \La^1 + \dots + \La^k$. 
Set $\bs \La^0 = (\La_1^0, \dots , \La^0_k)$ where 
\bea
\La_p^0\  = \ \La^p\ -\ \al(\bs l^p)
\eea
for $p = 1, \dots , k$.
Set $\La^0 = \La^0_1 + \dots + \La^0_k$.

Consider the tensor products
\bea
V_{\bs \La^0} & =& V_{\La^0_1} \otimes \cdots \otimes 
V_{\La^0_k} , 
\\
V_{\bs \La^p} & =& V_{\La^p_1} \otimes \cdots \otimes V_{\La^p_{n_p}} ,
\qquad
{\rm for}
\qquad
p = 1, \dots , k ,
\\
V_{\bs \La} &  = & V_{\bs \La^1} \otimes \cdots \otimes V_{\bs \La^k}  
\\
& =&
V_{\La^1_1} \otimes \cdots \otimes V_{\La^1_{n_1}} \otimes \cdots
\otimes V_{\La^k_1} \otimes \cdots \otimes V_{\La^k_{n_k}} . 
\eea
Let $S^0 $ be the tensor Shapovalov form on
$V_{\bs \La^0}$, $S^p$  the tensor Shapovalov form on 
$V_{\bs \La^p}$,  $S = S^1\otimes \cdots \otimes S^k$ the 
tensor Shapovalov form on $\V$.

To $p = 1, \dots , k$ and
$I = (i_1^1, \dots , i^1_{j_1};\ \dots ;\ i^{n_p}_1, \dots , i^{n_p}_{j_{n_p}})
\ \in P(\bs l^p, n_p)$ we associate a vector 
\bea
F_Iv_{\bs \La^p}\ =\ F_{i_1^1} \dots F_{i_{j_1}^1}v_{\La^p_1} \otimes \cdots
\otimes F_{i_1^{n_p}} \dots F_{i_{j_{n_p}}^{n_p}}v_{\La_{n_p}^p}
\eea
in $V_{\bs \La^p}[\La^p - \al(\bs l^p)]$. Assume that for $p = 1, \dots , k$
a singular vector
\bea
w_{\bs \La^p}\ = \ \sum_{I \in P(\bs l^p, n_p)}\
a^p_{I}\ F_Iv_{\bs \La^p} \ {} \ \in \ {} \s V_{\bs \La^p}[\La^p - \al(\bs l^p)]
\eea
is chosen. Here $a^p_{I}$ are some complex numbers.

To every $I = (i_1^1, \dots , i^1_{j_1};\ \dots ;\ i^{k}_1, \dots , i^{k}_{j_{k}}) \
 \in P(\bs l^0, k)$ we associate a vector 
\bea
F_I v_{\bs \La^0}\  =\ F_{i_1^1} \dots F_{i_{j_1}^1} v_{\La^0_1} 
\otimes \cdots
\otimes F_{i_1^k} \dots F_{i_{j_k}^k} v_{\La^0_k}
\eea
in $V_{\bs \La^0}[\La - \sum_{p=0}^k \al(\bs l^p)]$. Assume that
a singular vector
\bea
w_{\bs \La^0} \ =\ \sum_{I \in P(\bs l^0, k)}\
a^0_{I}\ F_I v_{\bs \La^0} \ {} \ \in \ {} \s V_{\bs \La^0}[\La - \sum_{p=0}^k
\al(\bs l^p)]
\eea
is chosen. Here $a^0_{I}$ are some complex numbers.

To every 
$I \in P(\bs l^0, k)$ we also associate a vector 
\bea
F_I w \  =\ F_{i_1^1} \dots F_{i_{j_1}^1} w_{\bs \La^1}\ 
\otimes \ \cdots\
\otimes\ F_{i_1^{k}} \dots F_{i_{j_k}^{k}} w_{\bs \La^k }
\eea
in $V_{\bs \La}[\La - \sum_{p=0}^k \al(\bs l^p)]$. Here
$F_{i_1^p} \dots F_{i_{j_p}^p} w_{\bs\La^p}$
denotes the action of
$F_{i_1^p} \dots F_{i_{j_p}^p}$ on the vector
$ w_{\bs\La^p}$ in the $\g$-module
$V_{\bs \La^p}$.

The vector
\bean\label{iterate}
\bs w\ =\  \sum_{I \in P(\bs l^0, k)}\
a^0_{I}\ F_I w \ {} \ \in \ {} \
\V[\La - \sum_{p=0}^k \al( \bs l^p)]
\eean
is called {\it the iterated singular vector with respect to the
singular vectors}
$ w_{\bs \La^0}, w_{\bs \La^1}, \dots , w_{\bs \La^k}$.
It is easy to see that $\bs w$ is a singular vector in $\V$.

\begin{lemma}\label{shap norm}
We have
\bea
S(\bs w, \bs w)\ =  \ \prod_{p=0}^k
S^p(w_{\bs \La ^p}, w_{\bs \La ^p})\ .
\phantom{aaaaaaaaaaaaaaaaaaaaaaaaa}
\square
\eea
\end{lemma}

\section{Asymptotics of Master Functions and Bethe Vectors}\label{ASYMP}

\subsection{Asymptotics of master functions}\label{sec as of master}
In this section we consider a master function
$\Phi (t,z, \bs \La, \bs l)$ and assume that parameters
$\bs \La, \bs l$ do not change while $z$
depends on a complex parameter $\epsilon$.
 We assume that $z$ has a limit as 
$\epsilon$ tends to zero. We study the limit of the master function, 
its critical points, and its Bethe vectors as $\epsilon$ tends to zero.

We use notations of Section \ref{iterated}.

Let $z = ( z_1, \dots , z_n)$. For $s = 1, \dots , n$ we assign 
the weight $\La^{p}_{s - n_1 - \dots - n_{p-1}}$
to the coordinate $z_s$ if 
\bean\label{4}
n_1 + \dots + n_{p-1} < s \leq
n_1 + \dots + n_{p} \ .
\eean
With this assignment we consider 
the master function $\Phi(t, z, \bs \La, \bs l)$ with
$t=(t_1, \dots , t_l)$.

Introduce the dependence of $z = (z_1, \dots , z_n)$ on new variables
$\epsilon$ and $(y^p_i)$ as follows.
Let $ y^0 = (y^0_1, \dots , y^0_k)$. For $p = 1, \dots , k$, let
$y^p = (y^p_1, \dots , y^p_{n_p})$. Let
$ y = (y^p_i)$ where $p = 0, \dots , k$ and $i = 1, \dots , k$ if $p = 0$
and $i = 1, \dots , n_p$ if $p = 1, \dots , k$.
Set
\bean\label{3}
z_s(y, \epsilon)\  =\  y^0_{p}\  + \ 
\epsilon \ y^{p}_{s - n_1 - \dots - n_{p-1}} , 
\eean
if $s$ satisfies \Ref{4}.

If the variables $y$ are fixed and $\epsilon \to 0$, then
the  coordinate $z_s(y, \epsilon)$ in \Ref{3} tends to
$y^0_p$ and the ratio $(z_s(y, \epsilon) - y^0_p)/\epsilon$ has the limit 
$y^{p}_{s - n_1 - \dots - n_{p-1}}$.

Let $z = z(y, \epsilon)$ be the relation given by formula  \Ref{3}.

We rescale the variables $t$ of the master function
$\Phi(t, z(y, \epsilon), \bs \La, \bs l)$  as follows.
Introduce new variables $u = (u^j_i)$ where
$j = 0, 1, \dots , k$ and $i = 1, \dots , l^j$.
If 
\bea
l_1 + \dots + l_{j-1} < i \leq 
l_1 + \dots + l_{j-1} + l_j^0 ,
\eea
then we set
\bean\label{1}
t_i \ =\ u^0_{l^0_1 + \dots + l^0_{j-1} + i - (l_1 + \dots + l_{j-1})} \ .
\eean
If
\bea
l_1 + \dots + l_{j-1} + l^0_j + \dots + l^{p-1}_j < i \leq 
l_1 + \dots + l_{j-1} + l^0_j + \dots + l^{p}_j ,
\eea
then we set
\bean\label{2}
t_i \ {} = \ {} y^0_p\ {} +\ {}  \epsilon\
 u^p_{l^p_1 + \dots + l^p_{j-1} + i - 
(l_1 + \dots + l_{j-1} + l^0_j + \dots + l^{p-1}_j)} \  .
\eean
Let $t = t(u, \epsilon)$ be the relation given by formulas  \Ref{1} and \Ref {2}.
The relation $t = t(u, \epsilon)$, given by formulas  \Ref{1} and \Ref {2},
will be called {\it the rescaling of variables $t$ with respect to the parameters
$\bs l^0, \dots , \bs l^k$} or simply {\it the 
$(\bs l^0, \dots , \bs l^k)$-type rescaling}.

We study the asymptotics of the function 
$
\Phi( t(u, \epsilon),    z(y, \epsilon), \bs \La, \bs l)  
$ as $\epsilon$ tends to zero.

To describe the asymptotics we use the master functions
$\Phi (u^p, y^p, \bs \La^p, \bs l^p)$, $p = 0, \dots , k$.
Here $u^p = (u^p_1, \dots , u^p_{l^p})$ for $p = 0, \dots , k$;\
$y^0 = (y^0_1, \dots , y^0_{k})$;\
$y^p = (y^p_1, \dots , y^p_{n_p})$ \linebreak
for $p = 1, \dots , k$;\
$\bs \La^p = (\La^p_1, \dots , \La^p_{n_p})$ for $p = 0, \dots , k$;\
$\bs l^p = (l^p_1, \dots , l^p_{r})$ \linebreak
for $p = 0, \dots , k$.

\begin{lemma}\label{as of master}
Let all the parameters $\La^j_i,\ l^j_i$ be fixed.
Fix a compact subset $K\subset \C^l\times \C^n$ in the $(u, y)$-space such that
the $y^0_1, \dots , y^0_k$ coordinates of points in $K$ are distinct.
Assume that $\epsilon$ tends to 0. Then
\bea
\Phi( t(u, \epsilon),    z(y, \epsilon), \bs \La, \bs l)
\ =\ 
\epsilon^{N(\bs \La,\ \bs l^1, \dots , \bs l^k)}\
(1 + \mathcal {O} (\epsilon, u, y) )\ 
\prod_{p = 0}^k \
\Phi (u^p, y^p, \bs \La^p, \bs l^p) .
\eea
Here 
$N(\bs \La,\ \bs l^1, \dots , \bs l^k)$ is a suitable constant. 
The function $\mathcal {O} (\epsilon, u, y) $ 
is holomorphic in $\C\times\C^l\times\C^n$ in a neighborhood of the set
$\{0\}\times K$ and $\mathcal {O} (\epsilon, u, y) |_{\epsilon = 0} = 0$.
\hfill
$\square$
\end{lemma}

\subsection{Asymptotics of critical points}\label{sec as of critical}
We keep notations of Section \ref{sec as of master}.

Let $y^0(*) = (y^0_1(*), \dots , y^0_{k}(*))$ be a point in $\C^k$ 
with distinct coordinates.
Let $u^0(*) = (u^0_1(*), \dots , u^0_{l^0}(*))$ be a non-degenerate
critical point
of the master function 
\linebreak
 $\Phi ( . , y^0(*), \bs \La^0, \bs l^0)$.

For $p = 1, \dots , k$ let
$y^p(*) = (y^p_1(*), \dots , y^p_{n_p}(*))$ be a point in $\C^{n_p}$ 
with distinct coordinates.
Let $u^p(*) = (u^p_1(*), \dots , u^p_{l^p}(*))$ be a non-degenerate
critical point
of the master function  $\Phi ( . , y^p(*), \bs \La^p, \bs l^p)$.

\begin{lemma}\label{lem as crit}
There exist unique functions $u^p_i(\epsilon)$, where
 $p = 0, \dots , k$ and $i = 1, \dots , k$ if $p = 0$
and $i = 1, \dots , n_p$ if $p = 1, \dots , k$,
with the following properties:
\begin{enumerate}
\item[$\bullet$]
The functions  $u^p_i(\epsilon)$ are holomorphic functions defined
in a neighborhood of  $\epsilon = 0$ in $\C$.

\item[$\bullet$]
We have\ $u^p_i(0) = u^p_i(*)$\ for all $p,\ i$.

\item[$\bullet$]
For all non-zero $\epsilon$ in a neighorhood of  $\epsilon = 0$
in $\C$ the point $u(\epsilon) = (u^p_i(\epsilon))$ is a non-degenerate
critical point of the function
$\Phi( t( u , \epsilon),    z(y(*), \epsilon), \bs \La, \bs l)$ with respect to 
the variables $u = (u^p_i)$.
\end{enumerate}
\end{lemma}

The lemma follows from Lemma \ref{as of master} with the
help of the implicit function theorem.

Let $u(\epsilon)$ be as in Lemma \ref{lem as crit}. Then
for small non-zero $\epsilon$, the point
\linebreak
$t(\epsilon) = t(u(\epsilon), \epsilon) \ \in \C^l$ is a non-degenerate critical
point of the master function
\linebreak
$\Phi( . ,    z(y(*), \epsilon), \bs \La, \bs l)$.
This family of critical points $t(\epsilon)$ of 
$\Phi( . ,    z(y(*), \epsilon), \bs \La, \bs l)$ will be called
\linebreak
{\it the family of critical points 
associated with the $(\bs l^0, \dots , \bs l^k)$-type rescaling
and originated at the
critical points $u^0(*), \dots , u^k(*)$ of the master functions
$\Phi ( . , y^0(*), \bs \La^0, \bs l^0)$, \dots , 
$\Phi ( . , y^k(*), \bs \La^k, \bs l^k)$, respectively}.

\subsection{Asymptotics of Hessians}\label{sec as of hessians}
If $f$ is a function of $t_1, \dots , t_n$ and 
$t(*)$ $=(t_1(*), \dots ,$ $ t_n(*))$
is a point, then the determinant
\bea
\det_{i,j = 1, \dots , n} 
\frac{\partial^2 f}{\partial t_i\, \partial t_j}( t(*) )
\eea
is called {\it the Hessian of $f$ at $t(*)$ with respect to variables 
$t = (t_1, \dots , t_n) $}
and denoted by ${\rm Hess}_t \, f(t(*))$.

\begin{lemma}\label{lem as of Hess}
Let $t(\epsilon)$ be the family of non-degenerate critical points
of the master function 
$\Phi( . ,    z(y(*), \epsilon), \bs \La, \bs l)$
associated with the $(\bs l^0, \dots , \bs l^k)$-type rescaling
and originated at the
critical points $u^0(*), \dots , u^k(*)$ of the master functions
$\Phi ( . , y^0(*), \bs \La^0, \bs l^0)$, \dots , \linebreak
 $\Phi ( . , y^k(*), \bs \La^k, \bs l^k)$, respectively. Then
\bea
\lim_{\epsilon\to 0}\ \epsilon^{ 2 (l^1 + \dots + l^k)}\
{\rm Hess}_t\,\ {\rm log}\ \Phi( t(\epsilon) ,    z(y(*), \epsilon), \bs \La, \bs l)
\ = \ \phantom{aaaaaaaaaaaaaaaaaaa}
\\
\phantom{aaaaaaaaaaaaaaaaaaa}
\prod_{p=0}^k\
{\rm Hess}_{u^p}\, {\rm log}\ \Phi(u^p(*), y^p(*), \bs \La^p, \bs l^p) \ .
\phantom{aaa}
\square
\eea
\end{lemma}

\subsection{Asymptotics of Bethe vectors}\label{sec as of vectors}
Let $t(\epsilon)$ be the family of non-degenerate critical points
of the master function 
$\Phi( . ,    z(y(*), \epsilon), \bs \La, \bs l)$
associated with the $(\bs l^0, \dots , \bs l^k)$-type rescaling
and originated at the
critical points $u^0(*), \dots , u^k(*)$ of the master functions
$\Phi ( . , y^0(*), \bs \La^0, \bs l^0)$, \dots , 
$\Phi ( . , y^k(*), \bs \La^k, \bs l^k)$, respectively.

Let 
\bea
\omega (t(\epsilon), z(y(*), \epsilon))\ \in
\ {} 
\s \V[\La - \sum_{p=0}^k \al(\bs l^p)]
\eea
be the Bethe vector corresponding to the critical point $t(\epsilon)$
of $\Phi( . ,    z(y(*), \epsilon), \bs \La, \bs l)$.

For $p = 0, \dots , k$ let 
\bea
\omega(u^p(*), y^p(*))\ {} \in\ {}
V_{\bs \La^p}[\La^p - \al(\bs l^p)]
\eea
be the Bethe vector corresponding to the critical point $u^p(*)$
of $\Phi( . , y^p(*), \bs \La^p, \bs l^p)$.

Let 
\bea
\bs \omega_{\omega_{\bs \La^0}, \omega_{\bs \La^1}, \dots , \omega_{\bs \La^k}} 
\ {} \in \ {} \s \V[\La -\sum_{p=0}^k\al(\bs l^p)]
\eea
be the iterated singular vector with respect to singular vectors
$\omega_{\bs \La^0}, \omega_{\bs \La^1}, \dots , \omega_{\bs \La^k}$.

\begin{lemma}\label{thm as of vactors}
We have
\bea
\lim_{\epsilon\to 0}\ \epsilon^{ l^1 + \dots + l^k}\
\omega(t(\epsilon), z(y(*), \epsilon))\ = \
\bs \omega_{\omega_{\bs \La^0}, \omega_{\bs \La^1}, \dots , \omega_{\bs \La^k}} \ .
\phantom{aaa}
\square
\eea
\end{lemma}

The lemma easily follows from the formula for the universal weight 
function by repeated application of the identity
\bea
\frac 1 {(t_i-t_j)(t_j-t_k)}\ =
\frac 1 {(t_i-t_k)(t_j-t_k)}\ +
\frac 1 {(t_i-t_j)(t_i-t_k)}\ .
\eea

\subsection{Asymptotics of hamiltonians}\label{sec as of hamiltonians}
In this section we 
keep the notations and assumptions of Section \ref{sec as of vectors}.

For $s = 1, \dots , n$ let $K_s(z) : \V \to \V$ be the Gaudin hamiltonian
associated with the tensor product $\V$ and a point $z \in \C^n$.
Let $c_s(\epsilon)$ be the eigenvalue on the Bethe eigenvector
$\omega(t(\epsilon), z( y(*), \epsilon ))$ of the operator 
$K_s(z( y(*), \epsilon ))$.

For $i = 1, \dots , k$, let $K_i(y^0(*)) : V_{\bs \La^0} \to
V_{\bs \La^0}$ be the Gaudin hamiltonian associated with the tensor product
$V_{\bs \La^0}$ and the point $y^0(*) \in \C^k$.
Let $c_i^0(u^0(*), y^0(*))$ be the eigenvalue on the Bethe eigenvector
$\omega(u^0(*), y^0(*))$ of the operator $K_i(y^0(*))$.

For $p = 1, \dots , k$ and $i = 1, \dots , n_p$, 
let $K_i(y^p(*)) : V_{\bs \La^p} \to
V_{\bs \La^p}$ be the Gaudin hamiltonian associated with the tensor product
$V_{\bs \La^p}$ and the point $y^p(*) \in \C^{n_p}$.
Let $c_i^p(u^p(*), y^p(*))$ be the eigenvalue on the Bethe eigenvector
$\omega(u^p(*), y^p(*))$ of the operator $K_i(y^p(*))$.

Consider the tensor product $\V$ as the tensor product
$V_{\bs \La^1} \otimes \dots \otimes V_{\bs \La^k}$
of $k$ $\g$-modules. 
For $i = 1, \dots , k$, consider the Gaudin hamiltonian
$\widehat K_i(y^0(*)) : \V \to \V$, 
\bea
\widehat K_i(y^0(*))\ = \sum_{j = 1, \ j\neq i}^k\
\frac {\Omega^{(i, j)}}{ y^0_i(*) - y^0_j(*)} \ ,
\eea
associated with those $k$ $\g$-modules and the point $y^0(*) \in \C^k$.
For $p = 1, \dots , k$ and $i = 1, \dots , n_p$, denote by
$\widehat K_i(y^p(*))^{(p)}$ the linear operator on 
$\V = V_{\bs \La^1} \otimes \dots \otimes V_{\bs \La^k}$
acting as  $K_i(y^p(*))$ on the factor $V_{\bs \La^p}$ and as the 
identity on other factors of that tensor product.

\begin{lemma} \label{lem as hamiltonians}
Let $s \in \{1, \dots , n\}$ and $s$ satisfies \Ref{4}.
If \ $n_p = 1$, then
\bea
\lim_{\epsilon \to 0}\ K_s(z(y^0(*), \epsilon))\ =\
\widehat K_p(y^0(*)) \ 
\eea
and
\bea
\lim_{\epsilon \to 0}\ c_i(\epsilon)\ =\
c_p^0(u^0(*), y^0(*)) \ .
\eea
If \ $n_p > 1$, then
\bea
\lim_{\epsilon \to 0}\ \epsilon\ K_s(z(y^0(*), \epsilon))\ =\
\widehat K_{i - ( n_1 + \dots + n_{p-1})} (y^p(*))^{(p)} \ 
\eea
and
\bea
\lim_{\epsilon \to 0}\ \epsilon\ c_i(\epsilon)\ =\
c^p_{i - ( n_1 + \dots + n_{p-1})} (u^p(*), y^p(*)) \ .
\phantom{aaa}
\square
\eea
\end{lemma}

\section{Norms of Bethe Vectors and Hessians}

\subsection{The $z$-dependence of the norm of a Bethe vector}
We use notations of Section \ref{master sec}.

 Fix a collection of weights
$\bs \La = (\La_1, \dots , \La_n)$ and a collection of non-negative
integers $\bs l = (l_1, \dots , l_r)$. Consider the master function
$\Phi(t, z, \bs \La, \bs l)$.

Let $z^0 = (z_1^0, \dots , z_n^0)$ be a point with distinct coordinates. 
Let $t^0 =   (t_1^0, \dots , t_l^0)$ be a non-degenerate critical point
of the master function $\Phi( . , z^0, \bs \La, \bs l)$.
By the implicit function theorem there exists a unique holomorphic
$\C^l$-valued function $t = t(z)$, defined in the neighborhood of 
$z^0$ in $\C^n$, such that $t(z)$ is a non-degenerate critical point
of the master function $\Phi( . , z, \bs \La, \bs l)$ and $t(z^0)=t^0$.
Let $\omega (t(z),z) \in \s \V[\Ll]$ be the corresponding Bethe vector. 
Let $S$ be the tensor Shapovalov form on $\V$.

\begin{theorem}[\cite{RV}] \label{const}
We have
\bean\label{const relation}
S ( \omega (t(z),z), \omega (t(z),z) )\ =\
C \ {\rm Hess}_t\ {\rm log}\ \Phi (t(z), z, \bs \La, \bs l) \ ,
\eean
where $C$ does not depend on $z$.
\end{theorem}

\begin{conjecture}[\cite{RV}]
\label{const conjecture} The constant $C$
in \Ref{const relation} is equal to $1$.
\end{conjecture}

The conjecture is proved for $\g = sl_2$ in \cite{V2}.
We prove the conjecture for $\g = sl_{r+1}$ in Theorem \ref{sl3 thm}.

\subsection{The upper bound estimate for the number of critical points}
Fix a collection of integral dominant $\g$-weights
$\bs \La = (\La_1, \dots , \La_n)$ and a collection of non-negative
integers $\bs l = (l_1, \dots , l_r)$. Consider the master function
$\Phi(t, z, \bs \La, \bs l)$ and its critical points with respect to $t$.
Recall the the group
$\bs\Sigma_{\bs l} = \Sigma_{l_1}\times \dots \times \Sigma_{l_r}$ 
acts on the critical set of $\Phi$.

\begin{theorem}\label{upper bound thm}
If $\La - \al(\bs l) $ is not a dominant integral $\g$-weight, then
the master function $\Phi( . , z, \bs \La, \bs l)$ does not have
isolated critical points, see Corollary 5.3 in \cite{MV5}.

If $\La - \al(\bs l) $ is  a dominant integral $\g$-weight, then
the master function $\Phi( . , z, \bs \La, \bs l)$ has only
isolated critical points, see Lemma 2.1 in \cite{MV2}.

If $\g = sl_{r+1}$ and
$\La - \al(\bs l) $ is a dominant integral $sl_{r+1}$-weight, then
the number of the $\bs\Sigma_{\bs l}$-orbits 
of critical points of the
master function $\Phi( . , z, \bs \La, \bs l)$, counted with multiplicities, is not
greater than the multiplicity of the irreducible
$sl_{r+1}$-module $V_{\La-\al(\bs l)}$
in the tensor product $\V$, see Theorem 5.13 
in \cite{MV2}.

If $\g = sl_{2}$, the weight
$\La - \al(\bs l) $ is a dominant integral $sl_{2}$-weight, 
and coordinates of the point $z = (z_1, \dots , z_n)$ are generic,
then the number of critical points of the
master function $\Phi( . , z, \bs \La, \bs l)$ is equal to
the multiplicity of the irreducible
$sl_{2}$-module $V_{\La-\al(\bs l)}$
in the tensor product $\V$. Moreover, in that case
all critical points are non-degenerate, 
see Theorem 1 in \cite{ScV}.
\end{theorem}

\subsection{Tensor products of two $sl_{r+1}$-modules if one of them
 is fundamental}
Let $\la$ be an integral dominant $sl_{r+1}$-weight, 
\ $w_1, \dots , w_r$ the fundamental $sl_{r+1}$-weights.
Set $e_1 = w_1, e_2 = w_2 - w_1, \dots , e_r = w_r - w_{r-1}, e_{r+1} = -w_r$.
For $p = 1, \dots , r$ we have
\bean\label{decomp formula}
V_{\la} \otimes V_{w_p} \ = \ \oplus_\mu\ V_\mu
\eean
where the sum is over all dominant integral weights $\mu$ such that
$\mu = \la + e_{i_1} + \dots + e_{i_p}$,\
$1\leq i_1 < \dots < i_p \leq r+1$.

For example if
$\la, \mu $ are dominant integral $sl_{r+1}$-weights,
then $V_\mu$ enters $ V_{\la} \otimes V_{w_1}$ if and only if
$\la = \mu - w_1 + \sum_{j=1}^i \al_j$ for some $i \leq r$.





Notice also that if $\la, \mu $ are dominant integral $sl_{r+1}$-weights,
then $V_\mu$ enters $ V_{\la} \otimes V_{w_r}$ if and only if
$\la = \mu - w_r + \sum_{j=i}^r \al_j$ for some $i \leq r$.


\bigskip

Consider the pair $\La = (\La_1, \La_2)$ where $\La_1$ is an
integral dominant $sl_{r+1}$-weight, and  $\La_2 = w_1$.
Write $\La_1 = \sum_{j=1}^r \la_jw_j$
for suitable non-negative integers $\la_j$.
Let $\bs l = (l_1, \dots , l_r) = 
(1, \dots , 1_i, 0_{i+1}, \dots , 0)$ for some $i \leq r$.
Assume that $\mu = \La_1 + w_1 - \al(\bs l)$ is an integral dominant weight.
Let $z^0 = (0, 1)$, and $t = (t_1, \dots , t_i)$.
Consider the master function $\Phi(t, z^0, \bs \La, \bs l)$.

Let $S$ be the tensor Shapovalov form on 
$V_{\La_1} \otimes V_{w_1}$.

\begin{theorem}[\cite{MV1}]\label{w_1 solution}
Under the above assumptions the function
$\Phi( . , z^0, \bs \La, \bs l)$ has exactly one critical point,
denoted by
$t^0 = (t^0_1, \dots , t^0_i)$. The critical point $t^0$
is non-degenerate. The coordinates of $t^0$ are given by the formula
\bean\label{coordinates}
t^0_j\ =\ \prod_{m=1}^j\
\frac
{\la_m + \dots + \la_i + i - m }
{\la_m + \dots + \la_i + i - m + 1} \ , 
\qquad
j = 1, \dots , i\ .
\eean
The Bethe vector $ \omega (t^0, z^0) \ \in 
\s V_{\La_1} \otimes V_{w_1} [ \La_1 + w_1 - \al(\bs l)]$,
corresponding to the critical point $t^0$, has the property
\bea\label{example}
S( \omega (t^0, z^0),  \omega (t^0, z^0))\ =\ 
{\rm Hess}_t\ {\rm log}\ \Phi (t^0, z^0, \bs \La, \bs l) \ .
\eea
\end{theorem}

Similarly 
consider the pair $\La = (\La_1, \La_2)$ where $\La_1$ is an
integral dominant $sl_{r+1}$-weight, and  $\La_2 = w_r$.
Let $\bs l = (l_1, \dots , l_r) = 
(0, \dots , 0_{i}, 1_{i+1}, \dots , 1)$ for some $i < r$. 
Assume that $\mu = \La_1 + w_r - \al(\bs l)$ is an integral dominant weight.
Let $z^0 = (0, 1)$, and $t = (t_1, \dots , t_{r-i})$.
Consider the master function $\Phi(t, z^0, \bs \La, \bs l)$.

Let $S$ be the tensor Shapovalov form on the tensor product
$V_{\La_1} \otimes V_{w_r}$.

\begin{theorem}[\cite{MV1}]\label{w_2 solution}
Under the above assumptions the function
$\Phi( . , z^0, \bs \La, \bs l)$ has exactly one critical point
, denoted by $t^0$. The critical point $t^0$
is non-degenerate.
The Bethe vector $ \omega (t^0, z^0)$ \ $\in 
\s V_{\La_1} \otimes V_{w_r}$ $ [ \La_1 + w_r - \al(\bs l)]$,
corresponding to the critical point $t^0$, has the property
\bea
S( \omega (t^0, z^0),  \omega (t^0, z^0))\ =\ 
{\rm Hess}_t\ {\rm log}\ \Phi (t^0, z^0, \bs \La, \bs l) \ .
\eea
\end{theorem}

The formulas for coordinates of the critical point in Theorem 
\ref{w_2 solution} can be easily deduced from formula
\Ref{coordinates}.

\subsection{Tensor products of two $sl_{4}$-modules if one of them
is the second fundamental}

If $\la, \mu $ are dominant integral $sl_{4}$-weights,
then $V_\mu$ enters $ V_{\la} \otimes V_{w_2}$ if and only if
$\la = \mu - w_2 + \delta$ where $\delta = 0$ or $\delta$
is one of the following five weights:
\bean \label{delta list}
\al_2,\ {}\  \al_1 + \al_2,\  {}\ \al_2 + \al_3,\  {}\ 
\al_1 + \al_2 + \al_3,\  {}\ \al_1 + 2\al_2 + \al_3\ .
\eean
For each  $\delta$ in \Ref{delta list}, 
write $\delta = l_1\al_1  + l_2\al_2 + l_3\al_3$
for suitable non-negative integers $l_i$. Set $\bs l = (l_1, l_2, l_3)$, \
$l = l_1 + l_2 + l_3$,\ $\bs \La = (\la, w_2)$, $z^0 = (0, 1)$, 
$t = (t_1, \dots , t_l)$.

Consider the master function $\Phi(t, z^0, \bs \La, \bs l)$.

\begin{theorem}\label{second fund}
Let $\la, \mu $ be dominant integral $sl_{4}$-weights, such that
$\la = \mu - w_2 + \delta$ and $\delta$ is one of the weights
in \Ref{delta list}. Then the function $\Phi ( . , z^0, \bs \La, \bs l)$
has exactly one critical point $t^0$. The critical point $t^0$ is non-degenerate.
The Bethe vector $\omega(t^0, z^0) \in \s V_\la\otimes V_{w_2} [\mu]$,
corresponding to $t^0$,  is a non-zero vector.
\end{theorem}

\begin{proof} If $\delta$ is 
$\al_2,\  \al_1 + \al_2$, or $ \al_2 + \al_3$, then the teorem follows from
Theorems \ref{w_1 solution} and \ref{w_2 solution}.

If $\delta$ is $\al_1 + \al_2 + \al_3$ or
$\al_1 + 2\al_2 + \al_3$, then the theorem is proved by direct verification.
Namely, let $\la = \la_1w_1 + \la_2w_2 + \la_3w_3$.
If $\delta = \al_1 + \al_2 + \al_3$, then one can check that
$t^0 = (t^0_1, t^0_2, t^0_3)$, where  
\bea
t^0_1 = \frac{ \la_1 (\la_1 + \la_2 + \la_3 + 2)}
{(\la_1 + 1) (\la_1 + \la_2 + \la_3 + 3)}\ ,
\qquad
t^0_2 = \frac{\la_1 + \la_2 + \la_3 + 2}
{\la_1 + \la_2 + \la_3 + 3}\ ,
\eea
\bea
t^0_3 = 
\frac{ \la_3 (\la_1 + \la_2 + \la_3 + 2)}
{(\la_3 + 1) (\la_1 + \la_2 + \la_3 + 3)}\ .
\eea
If $\delta$ is $\al_1 + 2\al_2 + \al_3$, then one can check that
$t^0 = (t^0_1, t^0_2, t^0_3, t^0_4)$, where
\bea
t^0_1 = 
\frac{(\la_1 + \la_2 + 1) (\la_1 + \la_2 + \la_3 + 2)}
{(\la_1 + \la_2 + 2) (\la_1 + \la_2 + \la_3 + 3)} ,
\qquad
t^0_4 = 
\frac{(\la_2 + \la_3 + 1) (\la_1 + \la_2 + \la_3 + 2)}
{(\la_2 + \la_3 + 2) (\la_1 + \la_2 + \la_3 + 3)}\  ,
\eea
\bea
&&
t^0_2 + t^0_3 - 2 = 
\\
&& \phantom{aaa}
-
\frac{(\la_1 + 2\la_2 + \la_3 + 4)
(\la_1\la_3 +  2\la_1\la_2 + 2\la_2\la_3 +
2(\la_2)^2 + 2\la_1 + 6\la_2 +  2\la_3 + 4)}
{(\la_2+1) (\la_1 + \la_2  + 2) (\la_2 + \la_3  + 2)
(\la_1 + \la_2 + \la_3 + 3)}\ ,
\eea
\bea
t^0_2 \, t^0_3  \ =\ 
\frac
{\la_2 (\la_1 + \la_2  + 1) (\la_2 + \la_3  + 1)
(\la_1 + \la_2 + \la_3 + 2)}
{(\la_2+1) (\la_1 + \la_2  + 2) (\la_2 + \la_3  + 2)
(\la_1 + \la_2 + \la_3 + 3)}\ .
\eea
One easily verifies the statements of the
theorem using those formulas.
\end{proof}

\section{Critical Points of the
$sl_{r+1}$ Master Functions with Frist and Last Fundamental Weights}

Let $\bs \La = (\La_1, \dots , \La_n)$ be a collection of $sl_{r+1}$-weights,
each of which is either the first or last fundamental, i.e. $\La_i \in \{w_1, w_{r}\}$. 
Let $\bs l = (l_1, \dots , l_r)$ be a sequence of non-negative integers such that
$\La - \al(\bs l)$ is integral dominant, here
$\La = \La_1 + \dots + \La_n$ and
$\al(\bs l) = l_1\al_1 + \dots + l_r\al_r$.

Consider the master function 
$\Phi(t, z, \bs \La, \bs l)$ where $t = (t_1, \dots , t_l), \ l = l_1 + \dots + l_r,$
and $z = ( z_1, \dots , z_n )$. Recall that the group
$\bs\Sigma_{\bs l} = \Sigma_{l_1}\times \dots \times \Sigma_{l_r}$ 
acts on the critical set of $\Phi( . , z, \bs \La, \bs l)$.

\begin{theorem}\label{sl3 fund weights}
For generic $z$ the following statements hold:
\begin{enumerate}
\item[(i)]
The number of \ $\bs\Sigma_{\bs l}$-orbits of
critical points of 
$\Phi( . , z, \bs \La, \bs l)$ is equal to the multiplicity of
the $sl_{r+1}$-module $V_{\La - \al(\bs l)}$ in the tensor product 
$V_{\bs \La}$.
\item[(ii)]
All critical points of $\Phi( . , z, \bs \La, \bs l)$
are non-degenerate. 
\item[(iii)] For every critical point $t^0$, the corresponding Bethe vector
$\omega (t^0, z)$ has the property:
\bea
S( \omega (t^0, z),  \omega (t^0, z))\ =\ 
{\rm Hess}_t\ {\rm log}\ \Phi (t^0, z, \bs \La, \bs l) \ .
\eea
\item[(iv)] 
The  Bethe vectors, corresponding to orbits of critical points of
$\Phi( . , z, \bs \La, \bs l)$,  form a basis in $\s \V[\Ll]$.
\end{enumerate}
\end{theorem}

\begin{proof}
The proof is by induction on $n$. If $ n = 2$, then the theorem follows from 
Theorems \ref{w_1 solution} and \ref{w_2 solution}.

Assume that Theorem \ref{sl3 fund weights} is proved for all tensor products of $n-1$
representations each of which is either the first or last 
fundamental. We prove Theorem \ref{sl3 fund weights} for
the tensor product $\V$ of $n$ given representations
$V_{\La_1}, \dots , V_{\La_n}$, each of which is either the first or last 
fundamental,  and the given sequence $\bs l = (l_1, \dots , l_r)$.
We  will use the notations and results of Sections \ref{iterated} and \ref{ASYMP}.

We may assume that $\La_n = w_1$. We may obtain that by either
reordering $\La_1, \dots , \La_n$ or using the automorphism of $sl_{r+1}$
which sends $E_i, H_i, F_i, \al_i, w_i$ to $E_{r+1-i}, H_{r+1-i},
F_{r+1-i}, \linebreak
\al_{r+1-i}, w_{r+1-i}$, respectively.

Introduce $n_1, \dots , n_k,\  \bs \La^1, \dots , \bs \La^k $ 
(as in Section \ref{iterated}) 
using the following formulas. Set
$k = 2$, $n_1 = n-1$, $n_2 = 1$, \
$\bs \La^1 = 
(\La_1, \La_2, \dots , \La_{n-1})$,
$\bs \La^2 = 
(\La_n)$,\ 
$V_{\bs \La^1}  = V_{\La_1} \otimes \cdots \otimes V_{\La_{n-1}}$,
\ $V_{\bs \La^2}  =  V_{\La_n}$,\
$V_{\bs \La}   =  V_{\bs \La^1} \otimes  V_{\bs \La^2} =
V_{\La_1} \otimes \cdots \otimes V_{\La_{n-1}} \otimes V_{\La_n}$.

Consider  the set $M'$ of the $r+1$ integral weights $\La - w_1 - \al(\bs l)$,
$\La - w_1 - \al(\bs l)+ \al_1$, \dots , 
$\La - w_1 - \al(\bs l)+ \al_1 + \dots + \al_r$. Denote by $M$ the subset of
all $\mu \in M'$ which are dominant.

Denote by $\mult (\mu; \la_1, \dots, \la_m)$ the multiplicity of $V_\mu$ in 
$V_{\la_1}\otimes \dots \otimes V_{\la_m}$. We have
\bea
\mult(\La - \al(\bs l); \La_1, \dots , \La_n)\ =\
\sum_{\mu \in M} \mult (\mu; \La_1, \dots , \La_{n-1})\ .
\eea

To prove parts (i-ii) 
of the theorem we will introduce a dependence of $z$ on $\epsilon$ so that
$z_1, \dots , z_{n-1}$ tend to 0 as $\epsilon \to 0$
and $z_n$ tends to 1. Using results of Section
\ref{ASYMP} we will construct non-intersecting sets of $\bs\Sigma_{\bs l}$-orbits
of critical points of $\Phi$, depending on $\epsilon$,  
labeled by $\mu \in M$, and consisting of  
$\mult (\mu; \La_1, \dots , \La_{n-1})$ elements each. Together
with Theorem \ref{upper bound thm} it will prove parts (i-ii).

More precisely, 
introduce the dependence of $z = (z_1, \dots , z_n)$ on the
new variables
$\epsilon$ and $y = (y^p_i)= (y^0_1, y^0_2, y^1_1, \dots , y^1_{n-1})$ as follows.
Set
\bean\label{z-dependence}
z_s(y, \epsilon)  &=&   y^0_1\  + \ 
\epsilon \ y^1_s ,
\qquad
s = 1, \dots , n-1 ,
\\
z_n(y, \epsilon)  &=&   y^0_2\ . 
\notag
\eean
Let $z = z(y, \epsilon)$ be the relation given by formula  \Ref{z-dependence}.
Set $y^0 = (y^0_1, y^0_2)$ and $y^1 =  (y^1_1, \dots , y^1_{n-1})$.

Introduce $r+1$ types of rescaling of coordinates $t$, cf. Section 
\ref{sec as of master}.

{\bf Type $0$ rescaling.} Set $\bs l^0 =  (0, \dots , 0),
\  \bs l^1 =  (l_1, \dots , l_r)$. 
Introduce new variables $u = (u^1_1, \dots , u^1_{l})$,
\bean\label{type I}
t_i \ {} = \ {} y^0_1\ {} +\ {}  \epsilon\
 u^1_i , 
\qquad i = 1, \dots , l\ .
\eean
This relation $t = t(u, \epsilon)$
will be called {\it the type $0$ rescaling of variables $t$}.
Set $u^0 = \emptyset$, $u^1 = (u^1_1, \dots , u^1_l)$.

{\bf Type $m$ rescaling,  $m = 1, \dots , r$.}
Set $\bs l^0 =  (1, \dots , 1_m, 0, \dots , 0),\
\bs l^1 = (l_1-1, \dots , \linebreak
l_m-1, l_{m+1}, \dots , l_r)$. 
Introduce new variables $u = (u^0_1, \dots , u^0_m,\ u^1_1, \dots , u^1_{l-m})$,
\bean\label{type II}
&{}&
\\
t_i  &=&   u^0_j \ ,
\qquad 
\phantom{aaaaaaaa}
{\rm if} \ {}
i = l_1 + \dots + l_{j-1} + 1 {}  \
{} {\rm for}\ {}j = 1, \dots , m \ ,
\notag
\\
t_i  &=&   y^0_1\ {} + \ {}  \epsilon\
 u^1_{i-j}\ , 
\qquad {\rm if}\ {}  l_1 + \dots + l_{j-1} + 1 < i \leq
l_1 + \dots + l_{j} \  {} {\rm for}\ {} j = 1, \dots , m \ , 
\notag
\\
t_i  &=&   y^0_1\ {} + \ {}  \epsilon\
 u^1_{i-m} \ , 
\qquad {\rm if}\ {}  l_1 + \dots + l_{m} < i \ . 
\notag
\eean
This relation $t = t(u, \epsilon)$
will be called {\it the type $m$ rescaling of variables $t$}.
Set $u^0 = (u^0_1, \dots , u^0_m)$, $u^1 = (u^1_1, \dots , u^1_{l-m})$.

We study the asymptotics of the function 
$\Phi( t(u, \epsilon),    z(y, \epsilon), \bs \La, \bs l)$
 as $\epsilon$ tends to zero for each of the $r+1$ rescalings.

To describe the asymptotics we use the master functions
$\Phi (u^p, y^p, \bs \La^p, \bs l^p)$, $p = 0, 1$.
Here the collections $\bs \La^1 = (\La_1, \La_2, \dots , \La_{n-1})$,
$\bs l^0, \bs l^1$, the variables 
$u^p$ and $y^p$ were already defined for each of the $r+1$
rescalings.
The collection $\bs \La^0$ is defined as follows.
For the type $0$ rescaling we set $\bs \La^0 = 
(\La^1 - \al(\bs l^1),  \La_n)$.
For the type $m$ rescaling with
$m = 1, \dots , r$, we set $\bs \La^0 = 
(\La^1 - \al(\bs l^1) + \al_1 + \dots + \al_m,  \La_n)$.

The master functions corresponding to the type $m$ rescaling will be provided
with the corresponding index:  
$\Phi_{m} (u^p, y^p, \bs \La^p, \bs l^p)$, $p = 0, 1$.

Let $y^1(*) = (y^1_1(*), \dots , y^1_{n-1}(*))$ be a point with distinct coordinates
such that:
\begin{enumerate}
\item[$\bullet$] For $m = 0, 1, \dots , r$, if
 $\La - w_1 -\al(\bs l)+ \al_1 + \dots + \al_m$ is dominant, then the
 master function
$\Phi_{m} (u^1, y^1(*), \bs \La^1, \bs l^1)$ has
$\mult (\La - w_1 - \al(\bs l) + \al_1 + \dots + \al_m; 
\La_1, \dots , \La_{n-1})$
distinct orbits of non-degenerate critical points satisfying
parts (iii-iv) of Theorem \ref{sl3 fund weights}.

\end{enumerate}
Such $y^1(*)$ exists according to the induction assumptions.

Consider the type $m$ rescaling with $m = 1, \dots , r$. 
Put $y^0(*) = (0, 1)$.
By Theorem \ref{w_1 solution} the function
$\Phi_m ( . , y^0(*), \bs \La^0, \bs l^0)$ has one critical point.
Denote the critical point by  $u^0(*) = (u^0_1(*), \dots , u^0_m(*))$.

Choose $\mult (\La - w_1 - \al(\bs l) + \al_1 + \dots + \al_m; \La_1, \dots ,
 \La_{n-1})$ critical points of \linebreak
$\Phi_{p} ( . , y^1(*), \bs \La^1, \bs l^1)$  lying in different 
$\Sigma_{l_1-1}\times \dots \times \Sigma_{l_m-1}\times \Sigma_{l_{m+1}} \times
\dots \times \Sigma_{l_r}$-orbits. Denote those critical points
by $u^1(*_j)$, $j = 1, \dots , 
\mult (\La - w_1 - \al(\bs l) + \al_1 + \dots + \al_m; \La_1, \dots ,
 \La_{n-1})$.
Let $t(\epsilon, j, m) \in \C^l$ be the family of critical points
of 
$\Phi( . ,    z(y(*), \epsilon), \bs \La, \bs l)$
associated with type $m$ rescaling and originated at the
critical points $u^0(*), u^1(*_j)$ of the master functions
$\Phi_{m} ( . , y^0(*), \bs \La^0, \bs l^0)$,
$\Phi_{m} ( . , y^1(*), \bs \La^1, \bs l^1)$, respectively,
see Section \ref{sec as of critical}.

Consider the type $0$ rescaling. Put $y^0(*) = (0, 1)$.
The function
$\Phi_{0} ( u^0 , y^0(*), \bs \La^0, \bs l^0)$ does not depend on
$u^0$.

Choose $\mult (\La - w_1 - \al(\bs l); \La_1, \dots ,
 \La_{n-1})$ critical points of 
$\Phi_{0} ( . , y^1(*), \bs \La^1, \bs l^1)$  lying in different 
$\Sigma_{l_1}\times \dots \times \Sigma_{l_r}$-orbits. Denote the critical points
by $u^1(*_j)$, $j = 1, \dots , 
\mult (\La - w_1 - \al(\bs l); \La_1, \dots ,
 \La_{n-1})$.
Let $t(\epsilon, j, 0) \in \C^l$ be the family of critical points
of \linebreak
$\Phi( . ,    z(y(*), \epsilon), \bs \La, \bs l)$
associated with type $0$ rescaling and originated at the
critical point $u^1(*_j)$ of the master function
$\Phi_{0} ( . , y^1(*), \bs \La^1, \bs l^1)$, see
Section \ref{sec as of critical}.

All together we constructed 
$\mult(\La - \al(\bs l); \La_1, \dots , \La_n)$ families of critical points of
$\Phi( . ,    z(y(*), \epsilon), \bs \La, \bs l)$.

The constructed families are all different. Indeed, the families corresponding to the
same rescaling are different by construction. The families, corresponding to
different rescalings are different because they have different limits 
as $\epsilon$ tends to 0. Now Theorem \ref{upper bound thm} implies
part (i).

All constructed critical points are non-degenerate by Lemma 
\ref{lem as crit}. This proves part (ii). Part (iii) is a direct
corollary of the induction assumptions, Theorems \ref{const},
\ref{w_1 solution}, and 
Lemmas \ref{thm as of vactors},  \ref{lem as of Hess}.

Part (iv) is a direct corollary of the construction and Lemma 
\ref{thm as of vactors}.
\end{proof}

\bigskip

Let $\bs \La = (\La_1, \dots , \La_n)$ be a collection of $sl_{4}$-weights,
each of which is  fundamental, i.e. $\La_i \in \{w_1, w_2, w_3\}$. 
Let $\bs l = (l_1, l_2 , l_3)$ be a sequence of non-negative integers such that
$\La - \al(\bs l)$ is integral dominant, here
$\La = \La_1 + \dots + \La_n$ and
$\al(\bs l) = l_1\al_1 + l_2\al_2 + l_3\al_3$.

Consider the master function 
$\Phi(t, z, \bs \La, \bs l)$ where $t = (t_1, \dots , t_l), \ l = l_1 + l_2 + l_3,$
and $z = ( z_1, \dots , z_n )$. Recall that the group
$\bs\Sigma_{\bs l} = \Sigma_{l_1}\times \Sigma_{l_2} \times \Sigma_{l_3}$ 
acts on the critical set of $\Phi( . , z, \bs \La, \bs l)$.

\begin{theorem}\label{sl4 fund weights}
For generic $z$ the following statements hold:
\begin{enumerate}
\item[(i)]
The number of \ $\bs\Sigma_{\bs l}$-orbits of
critical points of 
$\Phi( . , z, \bs \La, \bs l)$ is equal to the multiplicity of
the $sl_{4}$-module $V_{\La - \al(\bs l)}$ in the tensor product 
$V_{\bs \La}$.
\item[(ii)]
All critical points of $\Phi( . , z, \bs \La, \bs l)$
are non-degenerate. 
\item[(iii)] 
The  Bethe vectors, corresponding to orbits of critical points of
$\Phi( . , z, \bs \La, \bs l)$, are non-zero vectors and
 form a basis in $\s \V[\Ll]$.
\end{enumerate}
\end{theorem}

The proof of this theorem is  parallel to the proof of Theorem
\ref{sl3 fund weights} and is based on Theorem {second fund}.

\section{Norms of Bethe Vectors in the $sl_{r+1}$ Gaudin Models}
Let $\bs \La^0 = (\La_1^0, \dots , \La_k^0)$ be a collection of 
$sl_{r+1}$ integral dominant
weights. Let $\bs l^0 = (l_1^0, \dots , l_r^0)$ 
be a sequence of non-negative integers such that
$\La^0 - \al(\bs l^0)$ is integral dominant. Here
$\La^0 = \La_1^0 + \dots + \La_n^0$ and
$\al(\bs l^0) = l_1^0 \al_1 + \dots + l_r^0 \al_r$.

Consider the master function 
$\Phi(u^0, y^0, \bs \La^0, \bs l^0)$ where 
$u^0 = (u^0_1, \dots , u^0_{l^0}), \ l^0 = l_1^0 + \dots + l_r^0$,
and $y^0 = ( y^0_1, \dots , y^0_k )$. 

\begin{theorem}\label{sl3 thm}
Let $y^0(*) \in \C^k$ be a point with distinct coordinates.
Let $u^0(*)$ be a non-degenerate critical point of
$\Phi( . , y^0(*), \bs \La^0, \bs l^0)$. Let
$\omega(u^0(*), y^0(*)) \in \s V_{\bs \La^0}[\La^0 - \al(\bs l^0)]$ 
be the corresponding Bethe vector. 
Let $S^0$ be the tensor Shapovalov form on $V_{\bs \La^0}$. Then
\bea
S^0(\, \omega (u^0(*), y^0(*)), \ \omega (u^0(*), y^0(*))\,)\ =\ 
{\rm Hess}_{u^0}\ {\rm log}\ \Phi (u^0(*), y^0(*), \bs \La^0, \bs l^0) \ .
\eea
\end{theorem}

\begin{corollary}\label{non-zero vector}
The Bethe vector $ \omega(u^0(*), y^0(*))$ is a non-zero vector.
\end{corollary}

\begin{proof} We deduce Theorem \ref{sl3 thm} from
Theorem \ref{sl3 fund weights} using results of Section
\ref{ASYMP}.

It is known that for each dominant integral $sl_{r+1}$-weight $\la$, 
the multiplicity of $V_\la$ in $V_{w_1}^{\otimes n}$ is positive for
a suitable $n$.

For each $p = 1, \dots , k$ fix $n_p$ such that the multiplicity of 
$V_{\La_p^0}$ in $V_{w_1}^{\otimes n_p}$ is positive.
Set $\La^p = (w_1, \dots , w_1)$ where $w_1$ is taken $n_p$ times.
Denote by $S^p$ the tensor product Shapovalov form on $V_{w_1}^{\otimes n_p}$.

We have $n_p w_1 - \La_p^0 = l^p_1\al_1 + \dots + l^p_r\al_r$ where
$\bs l^p = (l^p_1, \dots , l^p_r)$ is a sequence of non-zero integers.
Set $l^p = l^p_1 +  \dots + l^p_r$, \
$y^p = (y^p_1, \dots , y^p_{n_p})$, \ 
$u^p = (u^p_1, \dots , u^p_{l^p})$. 
Consider the master function $\Phi(u^p, y^p, \bs \La^p, \bs l^p)$.
That master function satisfies conditions of Theorem \ref{sl3 fund weights}.
Hence there exists a point $y^p(*) \in \C^{n_p}$ with distinct coordinates
and a non-degenerate critical point $u^p(*) \in \C^{l^p}$ of the function
$\Phi( . , y^p(*), \bs \La^p, \bs l^p)$ such that the Bethe vector
$\omega(u^p(*), y^p(*)) \in \s V_{w_1}^{\otimes n_p}[\La_p^0]$ 
satisfies the identity:
\bea
S^p(\, \omega (u^p(*), y^p(*)),  \omega (u^p(*), y^p(*))\,)\ =\ 
\text{Hess}_{u^p}\ \text{log}\ \Phi (u^p(*), y^p(*), \bs \La^p, \bs l^p) \ .
\eea

Set $n = n_1 + \dots + n_k$, \
$\bs l = \bs l^0 + \dots + \bs l^k = 
(l^0_1 + \dots + l^k_1, \dots , l^0_r + \dots + l^k_r)$,
$l =  l^0 + \dots +  l^k$.
Set $z = (z_i^p)$, where $ p = 1, \dots , k,$ 
$i = 1, \dots , n_p$. Set $\bs \La = (\La^p_i)$, where 
$  p = 1, \dots , k, \ 
i = 1, \dots , n_p$, and $\La^p_i = w_1$. 
Assign the weight $\La^p_i$ to the variable $z^p_i$ for every $p, i$.
Set $t = (t_1, \dots , t_l)$. Consider the master function
$\Phi(t, z, \bs \La, \bs l)$.

Introduce the dependence of variables $z$ on variables $u$, 
$\epsilon$ by the formula: $z^p_i \ = \ y^0_p\ +\ 
\epsilon y^p_i $ for all $p, i$.
Introduce the $(\bs l^0, \dots , \bs l^k)$-rescaling of variables $t$ by formulas
\Ref{1} and \Ref{2}. Let $t(\epsilon) \in \C^l$ be 
the family of critical points 
associated with this rescaling
and originated at the critical points $u^0(*), \dots , u^k(*)$ 
of the master functions
$\Phi ( . , y^0(*), \bs \La^0, \bs l^0)$, \dots , 
$\Phi ( . , y^k(*), \bs \La^k, \bs l^k)$, respectively, see
Section \ref{sec as of critical}.

Let $\omega (t(\epsilon), z(y(*), \epsilon)) \in \s V_{w_1}^{\otimes n}$ 
be the corresponding Bethe vector.  
Let $ S$ be the tensor Shapovalov form on  $V_{w_1}^{\otimes n}$.
By Theorem \ref{sl3 fund weights} we have
\bea
S( \omega (t(\epsilon), z(y(*), \epsilon)), \omega (t(\epsilon), z(y(*), \epsilon)))\
 =\ 
{\rm Hess}_t\ {\rm log}\ \Phi (\omega (t(\epsilon), z(y(*), \epsilon)),
\bs \La, \bs l) \ .
\eea
Now by Lemmas  \ref{lem as of Hess}, \ref{thm as of vactors}, and \ref{shap norm} 
we may conclude that
\bea
S^0( \omega (u^0(*), y^0(*)),  \omega (u^0(*), y^0(*)))\ =\ 
{\rm Hess}_{u^0}\ {\rm log}\ \Phi (u^0(*), y^0(*), \bs \La^0, \bs l^0) \ .
\eea
\end{proof}

Similarly to Theorem \ref{sl3 thm} one can prove

\begin{theorem}\label{non-degenerate}
Let $t^0(*)$ be a critical point of 
$\Phi ( . , y^0(*), \bs \La^0, \bs l^0)$. Let
$\omega(u^0(*), y^0(*)) \in \s V_{\bs \La^0}[\La^0 - \al(\bs l^0)]$
be the corresponding Bethe vector. Assume that the number
\bea
S^0(\, \omega (u^0(*), y^0(*)), \ \omega (u^0(*), y^0(*))\,)
\eea 
is not equal to zero.
Then $t^0(*)$ is a non-degenerate critical point.
\end{theorem}

\begin{corollary}
Let $t^0(*)$ be a critical point of 
$\Phi ( . , y^0(*), \bs \La^0, \bs l^0)$ such that 
the corresponding Bethe vector
$\omega(u^0(*), y^0(*)) \in \s V_{\bs \La^0}[\La^0 - \al(\bs l^0)]$ 
is not equal to zero and belongs to the real part of $V_{\bs \La^0}$. 
Then $t^0(*)$ is a non-degenerate critical point.
\end{corollary}

The corollary follows from Theorem \ref{non-degenerate} since the Shapovalov 
form is positive definite on the real part of $V_{\bs \La^0}$. 

\medskip
\noindent
{\bf Example, cf \cite{RV}.} 
Let $\g = sl_2, \ \bs \La^0 = (w_1, w_1, w_1), \ \bs l^0 = (1),\
y^0(*) = (1, \eta, \eta^2)$, where $\eta = e^{2\pi i/3}$. Consider the master function
$\Phi(t, y^0(*), \bs \La^0, \bs l^0) = ((t_1)^3 - 1)^{-1}$. The point
$t^0(*)= (0)$ is the only critical point of $\Phi$. The critical point is degenerate.
 The corresponding Bethe vector 
\bea
\omega(u^0(*), y^0(*)) = &-& F_1v_{w_1}\otimes v_{w_1}\otimes v_{w_1}
\\
&-& \eta^2 \,v_{w_1}\otimes F_1v_{w_1}\otimes v_{w_1}
- \eta \,v_{w_1}\otimes v_{w_1}\otimes F_1v_{w_1}\ \in \ V_{\bs \La^0}
\eea
is a non-zero vector and
$S^0(\omega(u^0(*), y^0(*)), \omega(u^0(*), y^0(*)))= 1 + \eta^4 + \eta^2 = 0$.

\section{Transversality of   some   Schubert  Cells in 
  $Gr(r+1, \C_d[x])$}

In this section we formulate a corollary of Theorem
\ref{sl3 fund weights}.

Let $\mc V$ be a complex vector space of dimension $d + 1$  and
\be
\mathcal F=\{ 0\subset F_1\subset F_2\subset\dots\subset
F_{d+1} = \mc V \}, \qquad \dim F_i=i, 
\ee
a full flag in $\mc V$.  Let $Gr(r+1, \mc V)$ be the Grassmannian
variety of all $r+1$ dimensional subspaces in $\mc V$.

Let $\bs a =(a_1, \dots , a_{r+1})$, 
$ d - r \geq a_1\geq a_2\geq \dots \geq a_{r+1}\geq 0$,
be a non-increasing sequence of  non-negative integers.
Define the  {\it Schubert cell} $G^0_{\bs a}(\mc F)$,
associated to the flag $\mathcal F$ and the sequence $\bs a$, as the set
\bea
 \{V\in Gr(r+1, \mc V)\ | &&
\dim (V \cap F_{d- r + i-a_i}) = i ,
\\
&& \dim (V \cap F_{d- r + i-a_i-1}) = i-1, 
\ \text{for}\ 
i = 1, \dots , r+1 \} .
\eea
The closure $G_{\bs a}(\mc F)$ of the Schubert cell is called {\it
the Schubert cycle}. For a fixed flag $\mc F$, the Schubert cells 
form a cell decomposition of the Grassmannian.
The codimension of $G^0_{\bs a}(\mc F)$
in the Grassmannian is $|\bs a| = a_1 + \dots + a_{r+1}$. The cell corresponding to
$\bs a = (0, \dots , 0)$ is open in the Grassmannian.

Let $\mc V = \C_d[x]$ be the space of polynomials
of degree not greater than $d$, dim $\mc V = d + 1$. For any  $z \in \C \cup \infty$,
define a full flag in $\C_d[x]$,
\be
\mathcal F(z)\ = \ \{ 0 \subset F_1(z) \subset F_2(z) \subset \dots \subset
F_{d+1}(z) \}\ .
\ee 
For $z \in \C $ and any $i$, let $F_i (z)$ be the subspace of all
polynomials divisible by $( x - z )^{d+1-i}$. For any $i$,
let $F_i(\infty) $ be the subspace of all polynomials of degree less than $i$.

Thus,  for any sequence $\bs a$ and any $z \in \C \cup \infty$, 
we have a Schubert cell $G^0_{\bs a}(\mc F(z))$ in the Grassmannian
$Gr (r+1, \C_d[x])$ of all $r+1$-dimensional subspaces of $\C_d[x]$.

Let $V \in Gr (r+1, \C_d[x])$.
For any $z \in \C \cup \infty$, let $\bs a(z)$ be such a unique
sequence  that $V$ belongs to the cell $G^0_{\bs a (z)}(\mc F(z))$.
We say that a point $z \in \C \cup \infty$ is {\it a ramification point } for $V$, 
if  $\bs a(z) \neq (0, \dots , 0)$.

\begin{lemma}[\cite{MV2}]\label{codim}
For a basis $u_1, \dots , u_{r+1}$  in $V$, let
\ 
$W(u_1, \dots , u_{r+1}) \ =\ c \ \prod_{s=1}^n(x \ - \ z_s)^{m_s}$, $ c\neq 0$,
be its Wronskian.  Then
\begin{enumerate}
\item[$\bullet$] The ramification points for $V$ are the points $z_1, \dots , z_n$ and possibly
$\infty$.
\item[$\bullet$]  $| \bs a (z_s)| = m_s$ for every $s$.
\item[$\bullet$]  $| \bs a (\infty)| = (r+1)(d - r) - \sum_{s=1}^n m_s $.

\end{enumerate}
\end{lemma}

\begin{cor}\label{sum of cod} (Pl\"ucker formula)
We have
\bean\label{ram cond} 
\sum _{s=1}^n | \bs a (z_s)| \ + \ | \bs a (\infty)|\ =\ 
\dim \ Gr ( r+1, \C_d[x]) \ .
\eean
\end{cor}

A point $z \in \C$ is called {\it a base point} for $V$ if $u(z) = 0$ for every $u \in V$.
If $z\in \C$ is not a base point, then $a_{r+1}(z) = 0$.

Assume that ramification conditions 
$\bs a (z_1), \dots ,  \bs a (z_n), \bs a (\infty)$
are fixed at $z_1, \dots , z_n, \infty$, respectively, so that 
\Ref{ram cond} is satisfied and $a_{r+1}(z_s)=0$ for $s = 1, \dots , n$. 
The intersection number of Schubert cycles
$G_{\bs a (z_1)}(\mc F(z_1)), \dots , G_{ \bs a (z_n)}(\mc F(z_n)), 
G_{\bs a (\infty)}(\mc F(\infty))$ in the Grassmannian
$Gr (r+1, \C_d[x])$ can be described as follows.

Define integral dominant $sl_{r+1}$ weights $\La_1, \dots , \La_n, \La_\infty$
by the conditions
\bea
(\La_s,\al_i) = a_{r+1-i}(z_s) - a_{r+2-i}(z_s), 
\qquad 
(\La_\infty,\al_i) = a_{i}(\infty) - a_{i+1}(\infty),
\eea
for $i = 1, \dots , r$. The ramification conditions can be recovered from
$\La_1, \dots , \La_n, \La_\infty$ by the formula:
\bean\label{ramific conditions}
{}
\\
a_i(z_s) = (\La_s, \al_1 + \dots + \al_{r+1-i}),
\qquad 
a_i(\infty) = d - r - l_1 - (\La_\infty, \al_1 + \dots + \al_{i-1}), 
\notag
\eean
where
$l_1 = (\sum_{s=1}^{n}\La_s - \La_\infty, w_1)$ and $w_1$ is the
first fundamental weight.

According to Schubert calculus the intersection number of Schubert cycles
\linebreak
$G_{\bs a (z_1)}(\mc F(z_1)), \dots , 
G_{ \bs a (z_n)}(\mc F(z_n)),\  G_{\bs a (\infty)}(\mc F(\infty))$
is equal to the multiplicity of $V_{\La_\infty}$ in 
$V_{\La_1}\otimes \dots \otimes V_{\La_n}$, see \cite{Fu}.

\bigskip

Let $\bs \La = (\La_1, \dots , \La_n)$ be a collection of $sl_{r+1}$-weights,
each of which is either the first or last fundamental, i.e. $\La_p \in \{w_1, w_r\}$
for all $p$. 
Let $\bs l = (l_1, \dots , l_r)$ be a sequence of non-negative integers such that
$\La - \al(\bs l)$ is integral dominant, here
$\La = \La_1 + \dots + \La_n$ and
$\al(\bs l) = l_1\al_1 + \dots + l_r\al_r$. 

Fix a big positive integer $d$.

Let $z = (z_1, \dots , z_n)$ be a point in $\C^n$ with distinct coordinates.
By formula \Ref{ramific conditions} define ramification conditions
$\bs a (z_1), \dots ,  \bs a (z_n), \bs a (\infty)$ using the weights
$\La_1, \dots , \La_n, \La - \al(\bs l)$, respectively.
Thus $\bs a (z_s) = (1, \dots , 1, 0)$, if $\La_s = w_1$, 
$\bs a (z_s) = (1, 0, \dots , 0)$, if $\La_s = w_r$,
\linebreak
$\bs a (\infty) = 
(d - r - l_1, \ d - r + l_1 - l_2 - k_1,\ \dots ,\
d - r + l_{r-1} - l_r - k_1,\
d - r + l_{r} - k_1 - k_2)$, if $\La = k_1 w_1 + k_r w_r$.

\begin{theorem}\label{transver}
Under the above conditions on $\bs \La$,
for generic $z$ the intersection 
of Schubert cycles
$G_{\bs a (z_1)}(\mc F(z_1))$, $\dots ,$
$ G_{ \bs a (z_n)}(\mc F(z_n)),$
$ G_{\bs a (\infty)}(\mc F(\infty))$
in the Grassmannian $Gr(r+1, \C_d[x])$ consists of
$\mult (\La-\al(\bs l); \La_1, \dots , \La_n)$ distinct  points.
\end{theorem}

\begin{proof} By Theorem \ref{sl3 fund weights}
for generic $z$  the master function
$\Phi(t, z,  {\bs \La}, \bs l)$ has
$\mult (\La-\al(\bs l); \La_1, \dots , \La_n)$ distinct orbits
of critical points. According to Corollary 5.11 and 
Theorem 5.12 in \cite{MV2}
every orbit of critical points defines an intersection point
of Schubert cycles
$G_{\bs a (z_1)}(\mc F(z_1)), \dots , G_{ \bs a (z_n)}(\mc F(z_n)),
 G_{\bs a (\infty)}(\mc F(\infty))$
so that different orbits define different intersection points.
This proves the theorem since the intersection number of the cycles is 
equal to $\mult (\La-\al(\bs l); \La_1, \dots , \La_n)$.
\end{proof}

Note that the transversality properties of Schubert cycles
$G_{\bs a (z)}(\mc F(z))$  in the Grassmannian $Gr ( 2, \C_d[x])$
for arbitrary ramification conditions $\bs a(z)$ follow from the main theorem
in \cite{ScV}.

\end{document}